\documentclass{elsart}

\usepackage{amsmath,amsfonts,amssymb,graphicx,graphics}

\def\CT{\mathop{\mathrm{CT}}}

\def\xx{x^{-1}}
\def\yy{y^{-1}}
\def\ct{\mathop{\mathrm{CT}}}
\def\C{\mathcal{C}}
\def\G{\mathcal{G}}
\def\E{\mathcal{E}}
\def\H{\mathcal{H}}
\def\deg{\mathrm{deg}}

\def\figurename{Table}

\begin{document}

\begin{frontmatter}



\title{Enumeration of bilaterally symmetric $3$-noncrossing partitions}


\author[a]{Guoce Xin},
\ead{gxin@nankai.edu.cn}
\author{Terence Y. J. Zhang}
\ead{zyjcomb@gmail.com}
 \corauth[a]{Corresponding author.}
\address{
Center for Combinatorics, LPMC-TJKLC \\
Nankai University, Tianjin 300071, P.R. China}

\begin{abstract}
Sch\"{u}tzenberger's theorem for the ordinary RSK correspondence
naturally extends to Chen et. al's correspondence for matchings and
partitions. Thus the counting of bilaterally symmetric
$k$-noncrossing partitions naturally arises as an analogue for
involutions. In obtaining the analogous result for 3-noncrossing
partitions, we use a different technique to develop a {\sc Maple}
package for 2-dimensional vacillating lattice walk  enumeration
problems. The package    also   applies  to the hesitating case. As
applications, we find several interesting relations for some special
bilaterally symmetric partitions.
\end{abstract}

\end{frontmatter}

\section{Introduction}\label{sec:in}
A partition \(P\) of \([n]:=\{1,2,\dots,n\}\) is a collection of
nonempty subsets \(\{B_1, B_2, \ldots, B_k\}\), whose disjoint
union is \([n]\). The elements \(B_i\) are called \emph{blocks} of
\(P\).\(\) An important special class of partitions are (complete)
\emph{matchings} of \([2n]\), which are partitions of \([2n]\)
into
\(n\) two-element blocks.  
Every partition \(P\) of \([n]\) has a graph representation, called
\emph{partition graph}, obtained by identifying vertex \(i\) with
\((i,0)\) in the plane for \(i=1,\dots,n\), and drawing an arc
connecting \(i\) and \(j\) above the horizontal axis whenever \(i\)
and \(j\) are (numerically) consecutive in a block of \(P\). Such an
arc with \(i<j\) is called an edge \((i,j)\) of \(P\),  starting
from \(i\) and ending at \(j\). The vertices \(i\) and \(j\) are
called the \emph{left-hand endpoint} and the \emph{right-hand
endpoint} of the arc, respectively. A \emph{singleton} is the
element of a one-element block, and hence corresponds to an isolated
vertex in the graph. Conversely, a graph on the vertex set \([n]\)
is a partition graph if and only if each vertex is the left-hand
(resp., right-hand) endpoint of at most one edge. For a partition
\(P\) of \([n]\), let \(P^{refl}\) denote the partition obtained
from \(P\) by reflecting in the vertical line \(x=(n+1)/2 \).
Equivalently, \((i,j)\) is an arc of \(P\) if and only if
\((n+1-j,n+1-i)\) is an arc of \(P^{refl}\).

A sequence  \(\emptyset=\nu^0, \nu^1, \ldots, \nu^{2n}=\lambda \) of
Young diagrams is called a \emph{vacillating tableau} of shape
\(\lambda\) and length \(2n\) if (i) \(\nu^{2i+1}\) is obtained from
\(\nu^{2i}\) by doing nothing (i.e., \(\nu^{2i+1}=\nu^{2i}\)) or
deleting a square, and (ii) \(\nu^{2i}\) is obtained from
\(\nu^{2i-1}\) by doing nothing or adding a square.

In what follows, vacillating tableaux are always of shape
$\emptyset$ unless specified otherwise. Recently, Chen et al.
\cite{chen} established a bijection \(\phi\) from partitions   to
vacillating tableaux. Using their bijection, crossings and nestings
of a partition are characterized by its corresponding vacillating
tableau. When restricting to matchings, the image of \(\phi\)
becomes the set of oscillating tableaux. (see Appendix
\ref{sec:noncrossing} for definition).



 For a
vacillating tableau \(V\), reading \(V\) backward still gives a
vacillating tableau, denoted by \(V^{rev}\). Sch\"{u}tzenberger's
theorem for the ordinary RSK correspondence naturally extends to the
bijection $\phi$. The result for partitions is stated as follows.
\begin{thm}\label{sym-parti}
For any given partition \(P\) and  vacillating tableau \(V\),
\(\phi(P^{refl})=V^{rev}\) if and only if  \(\phi(P)=V\).
\end{thm}
This result and its analogy for matchings follows trivially from
Fomin's growth diagram language. See  \cite{Fomin-GRSKC}. The
matching case is due to Roby \cite{Roby} and the partition case is
due to Krattenthaler \cite{Krattenthaler}.

A vacillating tableau \(V\) is said to be \emph{palindromic} if \(V
= V^{rev}\). A partition \(P\) of \([n]\) is said to be
\emph{bilaterally symmetric} ({bi-symmetric for short}) if \(P =
P^{refl}\). Theorem \ref{sym-parti} implies that \(P\) is
bi-symmetric if and only if \(V(P)\) is palindromic.   The
enumeration of bi-symmetric partitions and matchings are not hard,
but turns out to be very difficult if we also consider the statistic
of \emph{crossing number} or \emph{nesting number}. A \(k\)-subset
\(\{(i_1,j_1),(i_2,j_2),\ldots,(i_k,j_k)\}\) of the edge set of a
partition \(P\) is said to be a \(k\)-crossing if
\(i_1<i_2<\cdots<i_k<j_1<j_2<\cdots<j_k\). A \(k\)-noncrossing
partition is a partition with no \(k\)-crossings. Some nice
properties on crossings and nestings of partitions and matchings
have been explored in \cite{chen}. Here we are interested in the
enumeration of these objects.

The number of \(k\)-noncrossing matchings was enumerated in
\cite{chen}, and the number of bi-symmetric \(k\)-noncrossing
matchings was enumerated in \cite{xin-weyl}. The number of
partitions is well-known to be the Bell number, but a formula for
the number of \(k\)-noncrossing partitions is only known for \(k=2\)
and \(k=3\). See \cite{Xin-3noncrossing}. The number of bi-symmetric
partitions was enumerated as the sequence A080107 in
 \cite{Sloane}. In this paper  we enumerate
bi-symmetric \(k\)-noncrossing partitions for \(k=2\) (In Appendix
\ref{sec:noncrossing}) and \(k=3\), which are the same as
palindromic vacillating tableaux of height bounded by \(k\) for
\(k=1\) and \(k=2\).

Let \(\widetilde{C}_3(n)\) be the number of bi-symmetric
\(3\)-noncrossing partitions of \([n]\). Then our main result is the
following.
\begin{prop}\label{bisym-main}
The numbers \(\widetilde{C}_3(2n)\) satisfy
\(\widetilde{C}_3(0)=1\), \(\widetilde{C}_3(2)=2\),
\(\widetilde{C}_3(4)=7\), and
\begin{multline}
27n(n+2)\widetilde{C}_3(2n)-3(7n^2+26n+27)\widetilde{C}_3(2n+2)-(7n^2+\\50n+84)\widetilde{C}_3(2n+4)
+(n+5)^2\widetilde{C}_3(2n+6)=0.\label{recfore}
\end{multline}
The numbers \(\widetilde{C}_3(2n+1)\) satisfy
\(\widetilde{C}_3(1)=1\), \(\widetilde{C}_3(3)=3\),
and
\begin{multline}
9(n^2+3n+2)\widetilde{C}_3(2n+1)-2(5n^2+30n+43)\widetilde{C}_3(2n+3)\\
+(n+4)(n+5)\widetilde{C}_3(2n+5)=0.\label{recforo}
\end{multline}
Equivalently, their associated generating functions
\(\G_e(t)=\sum_{n\geq0}\widetilde{C}_3(2n)t^n\) and
\(\G_o(t)=\sum_{n\geq0}\widetilde{C}_3(2n+1)t^n\) satisfy
\vspace{-3mm}
\begin{multline}-4-6t-6t^2+(4-12t-24t^2)\G_e(t)+(5t-29t^2-57t^3+81t^4){\frac {d}{d{t}}}\G_e(t)\\
+t^2(t-1)(3t+1)(9t-1){\frac {d^{2}}{d{t}^{2}}}\G_e(t)=0,
\label{deqfornce}
\end{multline}
\vspace{-13mm}
\begin{multline}
6+(-6+36t-18t^2)\G_o(t)+(-6t+50t^2-36t^3){\frac
{d}{d{t}}}\G_o(t)\\
-t^2(t-1)(9t-1){\frac
{d^{2}}{d{t}^{2}}}\G_o(t)=0.\label{deqfornco}
\end{multline}
\end{prop}

The above result is analogous to that for \(C_3(n)\), the number of
3-noncrossing partitions of \([n]\), in \cite{Xin-3noncrossing}. By
a similar way we represent the generating functions as certain
constant terms in Section \ref{sec:constant}. But the techniques
differs thereafter. In proving our result, we develop a {\sc Maple}
package in Section \ref{sec:maple} that applies to a class of two
dimensional vacillating lattice walk  enumeration problems. The
package is also extended to the hesitating case in Section
\ref{sec:hestating}.  As applications, we find several interesting
results for some special bi-symmetric partitions.

\section{Lattice Path Interpretations and Constant Term Expressions\label{sec:constant}}
In order to prove Proposition \ref{bisym-main}, we need to introduce
the lattice path interpretations. Let \(S\) be a subset of
\(\mathbb{Z}^k\). An \(S\)-\emph{vacillating lattice walk} of length
\(n\) is a sequence of lattice points \(p_0,p_1,\ldots,p_n\) in
\(S\) such that i) \(p_{2i+1}=p_{2i}\) or \(p_{2i+1}=p_{2i}-e_j\)
for some unit coordinate vector \(e_j\); ii) \(p_{2i}=p_{2i-1}\) or
\(p_{2i}=p_{2i-1}+e_j\) for some unit coordinate vector \(e_j\). We
are interested in two subsets of \(\mathbb{Z}^k\):
\(Q_k=\mathbb{N}^k\) of nonnegative integer lattice points and
\(W_k=\{(a_1,a_2,\ldots,a_k)\in \mathbb{Z}^k:a_1>a_2>\cdots>a_k\geq
0\}\) of Weyl lattice points. For two lattice points \(a\) and \(b\)
in \(W_k\) (or \(Q_k\)), denote by \(w_k(a,b,n)\) (or
\(q_k(a,b,n)\)) the number of \(W_k\) (or \(Q_k\))-vacillating
lattice walks of length \(n\) starting at \(a\) and ending at \(b\).
Let \(\delta=(k-1,k-2,\ldots,0)\).

Let \(C_k(n)\) be the number of \(k\)-noncrossing partitions of
\([n]\). The following consequence of Chen et. al's correspondence
\(\phi\) is the starting point of the enumeration for
\(3\)-noncrossing partitions, as well as for bi-symmetric
\(3\)-noncrossing partitions.
\begin{thm}[Chen et al., \cite{chen}]
The  number  \(C_{k+1}(n)\) equals \(w_k(\delta,\delta,2n)\), i.e.,
the number of closed \(W_k\)-vacillating lattice walks of length
\(2n\) from \(\delta\) to itself.
\end{thm}

By the correspondence \(\phi \), \(\widetilde{C}_3(n)\) is the same
as the number of palindromic vacillating tableaux of height  bounded
by \(2\) and length \(2n\), and is the same as the number of
palindromic \(W_2\)-vacillating lattice walks of length \(2n\) that
start and end at \((1,0)\). Since such walks are palindromic, it is
sufficient to consider only the first \(n\) steps of the lattice
walks. We have
\begin{align}\widetilde{C}_3(n)=\sum_{b\in W_2} w_2((1,0),b,n).
\label{latticC3}
\end{align}

Let us introduce the basic idea for solving the problem of
determining \(C_3(n)\), where the \(Q_2\)-vacillating lattice walks
 starting and ending at \((1,0)\) are considered. The same idea applies to
 determining \(\widetilde{C}_3(n)\).

It was shown in \cite{Xin-3noncrossing} by using the reflection
principle that \[w_k(a,b,n)=\sum_{\pi\in
\mathfrak{S}_k}(-1)^{\pi}q_k(\pi(a),b,n),\] where \((-1)^{\pi}\) is
the sign of \(\pi\) and
\(\pi(a_1,a_2,\ldots,a_k)=(a_{\pi(1)},a_{\pi(2)},\ldots,a_{\pi(k)})\).
Thus the enumeration of \(w_k(\delta,\delta,2n)\) reduces to that of
\(q_k(a,\delta,2n)\). Denote by \(a_{i,j}(n)=q_2((1,0),(i,j),n)\).
Let
\begin{align*}
F_e(x,y;t)&=\sum_{i,j,n\ge 0} a_{i,j}(2n)x^iy^jt^{2n}\\[-20pt]
\intertext{and}\\[-20pt]
 F_o(x,y;t)&=\sum_{i,j,n\ge 0} a_{i,j}(2n+1)x^iy^jt^{2n+1}
\end{align*}
be respectively the
 generating functions of lattice walks of even and odd length.
By a step by step construction, one can set up functional equations
for \(F_e(x,y;t)\) and \(F_o(x,y;t)\) and reduces the problem to
solving the following functional equation:
\begin{align*}
K(x,y;t^2) F_o(x,y;t)/t&= {x(1+\xx+ \yy)- \xx V_e(y;t^2)-\yy
H_e(x;t^2)} ,
\end{align*}
where
  \(V_e(y;t^2)\) and \(H_e(x;t^2)\)
are respectively the generating functions for lattice walks of even
length that start  at \((1,0)\) and end on the vertical and
horizontal axis, and   the \emph{kernel} of the equation
\(K(x,y;t)\)  is given by
\begin{align*} K(x,y;t)&={1-t(1+x+y)(1+\xx+\yy)}.
\end{align*} By the
{\em obstinate kernel method\/} of \cite{bousquet-motifs}, one can
finally obtain the generating function $\mathcal{C}(t)$ of $C_3(n)$
as
\begin{align}\label{e-Gtfirst}
\mathcal{C}(t)=\ct_x \;\big( (x^{-2}-x^2)(x^2+
(x^{-2}+x+x^2)Y+(x^{-3}-\xx) Y^2-x^{-2} Y^3)\big),
\end{align}
where the operator \(\ct_x\) extracts the constant term in \(x\) of
series in \(\mathbb{Q}[x,\xx][[t]]\) and  \(Y=Y(x;t)\) is the unique
power series in \(t\) satisfying \(Y=t(1+x+Y)(1+(1+x^{-1})Y)\) given
by
\begin{align}
Y&= \frac {1-(\xx+3+{x})t- \sqrt {\left(1-(1+x+ \xx)t\right)^2-4t}
}{ 2\left( 1+\xx \right) t}
=(1+x)t
+\cdots .\label{e-Y}
\end{align}  We shall mention that all this is done in the ring
\(Q[x,\xx,y,\yy][[t]]\) of formal power series in \(t\) with
coefficients Laurent polynomial in \(x\) and \(y\).

This idea works in a similar way for lattice walks starting from a
set of points and ending at \((1,0)\). For a set \(A\) of points, we
denote by \(A(x,y)=\sum_{(i,j)\in A} x^iy^j\) its generating
function. Let \(C^A_3(n)\) be the number of \(W_2\)-lattice walks of
length \(2n\) starting from points in \(A\) and ending at \((1,0)\),
and let \(\C^A(t)\) be the generating function of \(C^A_3(n)\). For
instance, \(A_1(x,y)=x\) corresponds to the point \((1,0)\) and
hence \(C^{A_1}_3(n)=C_3(n)\) and \(\C^{A_1}(t)=\C(t)\). For general
\(A\), with \(Y\) as in \eqref{e-Y} the result of \cite[Section
2.7]{Xin-3noncrossing} for \(\C^A(t)\) can be summarized as follows.
\begin{prop}\label{p-A}
For any set  \(A\) of lattice points in \(W_2\), we have
\begin{multline}
\mathcal{C}^A(t)=\CT_x\big((x^{-2}-x^2)( ( x+Y+xY ) A (x,Y) -
(x^{-1} Y+ Y
+x^{-1} Y^2) A (x^{-1} Y, Y)\\
 + (x^{-1} Y +x^{-1} +x^{-2} Y)
 A(x^{-1} Y, x^{-1}))\big).\label{ctfunction}
\end{multline}
\end{prop}

Now it is natural to let \(A_2(x,y)=\frac{x}{(1-x)(1-xy)}\), which
 corresponds to the set of all  points in \(W_2\). We shall also consider
the following two closely related cases: $A_3(x,y)=x/(1-x)$
corresponds
to the $x$-axis in \(W_2\); 
 $A_4(x,y)=x/(1-xy)$ corresponds to the diagonal in  \(W_2\).
Define \(e(n)=w_2((1,0),A_3,n)\) and \(h(n)=w_2((1,0),A_4,n)\). Then
at the same \(e(n)\) (resp., $h(n)$) is   the number of bi-symmetric
\(3\)-noncrossing partitions on \([n]\) whose central Young diagrams
consist of at most one row (resp., two rows of squares of equal
length including \(\emptyset\)).

Although our lattice walks  for \(\widetilde{C}_3(n)\) always start
from \((1,0)\), which is   different from that in Proposition
\ref{p-A}, we will still use the formulas for \(\C^A(t)\)  by means
of the following two observations:
\begin{enumerate}\item[1)]
\(w_2(a,b,2n)=w_2(b,a,2n)\), since \(W_2\)-vacillating lattice walks
of even length are still \(W_2\)-vacillating if read {backward}.
Thus by \eqref{latticC3} \(\widetilde{C}_3(2n)=C^{A_2}_3(2n)\), and
similarly \(e(2n)=C^{A_3}(2n)\) and   \(h(2n)=C^{A_4}(2n)\).

 \item[2)] By the step by step construction we have
 \[w_2(a,b,2n+1)=w_2(a,b,2n)+w_2(a,b+(1,0),2n)+w_2(a,b+(0,1),2n).\]
However, we must take care of the boundary cases. A careful study
yields
\begin{equation*}\sum_{b\in
W_2}w_2(a,b,2n+1)=w_2(a,(1,0),2n)+2\cdot \sum_{b\in
A_2'}w_2(a,b,2n)+3\cdot \sum_{b\in A_2''}w_2(a,b,2n),
\end{equation*}
where 
\(A_2'(x,y)=x^2/(1-x)+x^2y/(1-xy)\) and
\(A_2''(x,y)=x^3y/((1-x)(1-xy))\).

\begin{figure}[hbt]
\begin{center}
\begin{tabular}{ c|c }
 The ending set (of odd length ) & generating function set (of even length)\\
\hline
$A_2$ (the set of all  points in \(W_2\))&$A_1(x,y)+2A_2'(x,y)+3A_2''(x,y)$\\
$A_3$ (the $x$-axis  in \(W_2\)) &$(1+x+xy)A_3(x,y)$\\
$A_4$ (the diagonal  in \(W_2\)) &$(1+x)A_4(x,y)$\\
\end{tabular}
\caption{\label{tab1} Reducing the length from $2n+1$ to $2n$ by the
step by step construction.}
\end{center}
\end{figure}
\end{enumerate}


In summary, with   \(\G_e(t)\) and \(\G_o(t)\) as stated in
Proposition \ref{bisym-main}  we have
\begin{align*} \G_e(t) &=\C^{A_2(x,y)}(t),\\[-20pt]
\intertext{and}\\[-20pt]
\G_o(t)&=\C^{A_1(x,y)+2A_2'(x,y)+3A_2''(x,y)}(t).
\end{align*}
Then by Proposition \ref{p-A}, \(\G_e(t)\) and \(\G_o(t)\) can be
represented as certain constant terms. The cases for the other
\(A\)'s are similar. Such constant terms will be systematically
dealt with by the {\sc Maple} package developed in Section
\ref{sec:maple}.

Several  interesting results  can be obtained  for the $A_3$ and
$A_4$ cases similarly.
\begin{prop}\label{twice-relation}
For \(n\ge 1\), we have \(e(2n)=2 \cdot h(2n)\). Moreover
\(h(2)=1\), \(h(4)=3\)   and
\begin{equation} 9n(n+3)h(2n)-2(5n^2+26n+30)h(2n+2)+(n+4)(n+5)h(2n+4)=0.\end{equation}
\end{prop}
The proposition can be established by the following two differential
equations, which can be  easily shown by our package. \vspace{-3mm}
\begin{multline}-6-6t+(6-18t)\C^{A_3}(t)+(6t-42t^2+36t^3)\frac{d}{d
t}\C^{A_3}(t)\\+t^2(t-1)(9t-1)\frac{d^2}{d
t^2}\C^{A_3}(t)=0,\label{deqforxaxis}
\end{multline}
\vspace{-10mm}
\begin{multline} -6+6t+(6-18t)\C^{A_4}(t)+(6t-42t^2+36t^3)\frac{d}{d
t}\C^{A_4}(t)\\
+t^2(t-1)(9t-1)\frac{d^2}{d
t^2}\C^{A_4}(t)=0.\label{deqfordiag}
\end{multline}

\begin{prop}\label{vaci-odd}
For \(n\ge 0\), we have \(e(2n+1)=e(2n+2)/2\), and \(h(1)=1\),
\(h(3)=2\)   and
\begin{equation} 9(n+2)^2h(2n+1)-(10n^2+62n+93)h(2n+3)+(n+5)^2h(2n+5)=0.\label{oddA4}\end{equation}
\end{prop}

The sequence $(h(2n+1))_{n\geq 0}$ appears as A005802 in
\cite{Sloane}. This suggests that $h(2n+1)=u_{n+1}$, the number of
1234-avoiding permutations of length $n+1$.  It is  easy to check
that Equation \eqref{oddA4} coincides with    the formula given by
Mihailovs in the comments of A005802.
%

We conclude this subsection by some asymptotic estimates in Table
\ref{tab2}.

\medskip
\begin{figure}[hbt]
\begin{center}
\begin{tabular}{c|llllllllllll}
$n$&0&1&2&3&4&5&6&\(\rightarrow\)&\(\infty\)\\
\hline
\(w_2((1,0),A_1,2n)\)&1& 1& 2& 5& 15& 52& 202& \(\sim\)& \(\kappa_1\cdot9^n/n^7\)\\
\(w_2((1,0),A_2,2n)\)&1&2&7&30&148&806&4716&\(\sim\)&\(\kappa_2\cdot 9^n/n^3\)\\
\(w_2((1,0),A_2,2n+1)\)&1&3&12&57&303&1743&10629&\(\sim\)&\(\kappa_3 \cdot 9^n/n^3\)\\
\(w_2((1,0),A_3,2n)\)&1&2&6&22&94&450&2346&\(\sim\)&\(\kappa_4\cdot 9^n/n^4\)\\
\(w_2((1,0),A_4,2n)\)&1&1&3&11&47&225&1173&\(\sim\)&\(\kappa_5\cdot
9^n/n^4\)\\
\end{tabular}
 \caption{\label{tab2}
The first several numbers of vacillating
lattice walks and their asymptotic estimate, where
$\kappa_1\approx 1691.643,\kappa_2\approx 3.719,\kappa_3\approx
11.156$, and $\kappa_4=2\kappa_5\approx 16.732$.
 }
\end{center}
\end{figure}

\section{Determine the Constant Terms by a  Maple Package \label{sec:maple}}
In this section we will develop a  {\sc Maple} package to deal with
constant term expressions for \(\C^{A}(t) \). Our proof is based on
the idea of Lipshitz \cite{lipshitz-df}, but for our particular
problem we find a much smaller bound for the degree of the
D-finiteness. Moreover, this bound is for a large class of power
series and can be carried out by  {\sc Maple}. We find it better to
work in the filed \(\mathbb{Q}((x))((t))\) of iterated Laurent
series, which is also the field of Laurent series in \(t\) with
coefficients Laurent series in \(x\). See
\cite{xin-iterate,xin-residue}  for other applications of this
field.

Many objects are easy to describe using \[u=(\xx+2+x)=\xx(1+x)^2.
\] Let \[\Delta\equiv\Delta(x,t)=\sqrt {\left(1-(1+x+
x^{-1})t\right)^2-4t}=\sqrt {\left(1-(u-1)t\right)^2-4t}.\] Then it
is easy to see that
\(\mathbb{Q}(x,t,\Delta)=\mathbb{Q}(x,t)\oplus\mathbb{Q}(x,t)\Delta\).
Since
\begin{align}Y=\frac{1}{2}\,{\frac {x- \left( 1+{x}^{2}+3\,x \right) t-x\Delta}{t \left( 1+x
 \right) }},\label{Y}\end{align}
\(\C^A(t)\) can be written as \[\C^A(t)=\CT_xT_0+\CT_xT_1\Delta\]
for some \(T_0, T_1 \in \mathbb{Q}(x,t)\). In our study, the series
\(A(x,y)\) is always in the form of \(P(x,y)/\big((1-x)(1-xy)\big)\)
for some polynomial \(P(x,y)\). Consequently the rational functions
\(T_0\) and \(T_1\)   may have \(x,1+x,D_1,D_2,\) and \(D_3\) (but
no more) as denominators, where
\[D_1= -2t+(1-5t)x-2tx^2=x(1-(2u+1)t),\]
  \vspace{-7mm}
\[D_2= -t-2tx-(3t-1)x^2-2tx^3-tx^4=x^2(1-(u-1)^2 t),\]
%
\[D_3=
t^2+(2t^2-2t)x+(3t^2-6t+1)x^2+(2t^2-2t)x^3+t^2x^4
=x^2(\left(1-(u-1)t\right)^2-4t).\]

Of course one can write everything in terms of \(Y\), but using
\(\Delta\) may significantly simplify the proof because the
derivatives of \(\Delta\) have simple expressions. Notice that
\(D_3=\Delta^2x^2\), and  we have
\begin{align}
 \frac{\partial }{\partial x}
\Delta(x,t)&= \frac{t^2x^4+(t^2-t)x^3-(t^2-t)x-t^2}{xD_3}\Delta,
\label{pardiffDtox}
 \\
\frac{\partial }{\partial t} \Delta(x,t)&
=\frac{tx^4+(2t-1)x^3+(3t-3)x^2+(2t-1)x+t}{D_3}\Delta.
\label{pardiffDtot}
\end{align}

Let \(\mathcal{L}\) be the (finite) \(\mathbb{Q}(t)\)-linear span of
\[\left\{\CT_x \frac{L(x,t)\Delta}{(1+x)^{p}D_1^qD_2^rD_3^s}\mid
(p,q,r,s)\in \mathbb{Z}^4, L(x,t) \mathrm{\ is\ a\ Laurent\
polynomial\ in\ } x \right\}.\] We shall devote ourselves to prove
the following result.

\begin{prop}\label{laurent-prop}
The linear span \(\mathcal{L}\) is of dimension
 at most 3. More precisely, for any given \(L(x,t), p,q, r\)
and \(s\), there exists a procedure to find rational functions
\(R(t),P(t),Q(t)\in \mathbb{Q}(t)\) such that
\begin{equation}\CT_x
 \frac{L(x,t)\Delta}{(1+x)^{p}D_1^qD_2^rD_3^s} =R(t)+\CT_x(P(t)+Q(t)x)\Delta.\label{eqforprop}\end{equation}
\end{prop}
It is clear that \(\mathcal{L}\) is closed under taking derivatives
with respect to \(t\). Thus we have
\begin{cor}Every element in \(\mathcal{L}\) is \(D\)-finite of
 order at most 2.
 \end{cor}

The basic idea for proving Proposition \ref{laurent-prop} is to
use the well-known formula
\begin{align}\CT_x x\frac{\partial}{\partial
x}F(x,t)=0,  \text{ for all } F(x,t)\in
\mathbb{Q}((x))((t))\tag{P1}\label{P1}\end{align} to reduce
elements of \(\mathcal{L}\) into simple form.
 We need the following lemma.
\begin{lem}\label{eliminat-lem}
a) For all  $k\in \mathbb{Z}$, we have \vspace {-3mm}
\begin{align}
\CT_x(x^k-x^{-k})\Delta&=0,\tag{P2}\label{P2}\\
\CT_x   (x^k-x^{-k})\frac{x\Delta} {D_1}&=0, \tag{P3}\label{P3}\\
\CT_x   (x^k-x^{-k})\frac{x^2\Delta}{D_2}&=0. \tag{P4}\label{P4}
\end{align}
b) \vspace {-6mm}
\begin{align}
 \CT_x\frac{1-x}{1+x}\Delta&=1-t.\tag{P5}\label{P5}\\
\CT_x \frac{1-x}{1+x}  \frac{x\Delta}{D_1}&=1.\tag{P6}\label{P6}\\
\CT_x\frac{1-x}{1+x} \frac{x^2\Delta}{D_2}&=1.\tag{P7}\label{P7}
 \end{align}
c)\vspace {-6mm}
\begin{align}\CT_x
(1-x)(1-3xt)\Delta=(1-t)^2.\tag{P8}\label{P8}\end{align}
\end{lem}
{\bf{Proof.}} For brevity and similarity, we only prove
(\ref{P4},\ref{P7},\ref{P8}). Using the easy fact \[ \CT_x
F(x,t)=\CT_x F(\xx,t), \text{ if }  F(x,t)\in
\mbox{\(\mathbb{Q}\)}[x,1/x][[t]], \] we can prove \eqref{P4} by
letting \(F(x,t)=\Delta(x,t)/(D_2/x^2)\) and observing
\(F(x,t)=F(\xx,t)\).

For part b), we use Jacobi's change of variable formula
\cite{xin-residue} in the one variable case:
\begin{thm}[Jacobi's Residue Formula] Let \(y=f(x)\in \mathbb{C}((x)) \) be a Laurent
series and let \(b\) be the integer such that \(f(x)/x^b\) is a
formal power series with nonzero constant term. Then for any formal
series \(G(y)\) such that the composition \(G(f(x))\) is a Laurent
series, we have \vspace{-2mm}
\begin{align}
 \CT_x G(f(x)) \frac{x}{f}\frac{\partial
f}{\partial x} =b \CT_y G(y). \label{e-Jacobi}
\end{align}
\end{thm}

We make the change of variable by \(f(x)=u=\xx+2+x=\xx(1+x)^2\) with
\(b=-1\). It is worth mentioning that the \(y\) on the right-hand
side of \eqref{e-Jacobi} is understood the same as \(\xx\) (or very
large). For instance, $G(y)=1/(1-y)$ should be expanded as
\(1/(-y(1-1/y))=\sum_{n\ge 0}-y^{-1-n}\). See \cite{xin-residue} for
detailed explanation. Though this understanding is not used in our
calculation since \(G(y) \) will be taken as Laurent polynomials, it
is crucial if we make a more natural change of variable by
\(f(x)=\xx+1+x \).

Direct calculation shows that \[ \frac{x}{u}\frac{\partial
u}{\partial x}=\frac{x^2}{(1+x)^2} (1-x^{-2}) = -\frac{1-x}{1+x}.
\] Thus Jacobi's Residue Formula gives us the following equality
\[ \CT_x G(u(x)) \frac{1-x}{1+x} =\CT_u G(u).\]  Noticing that
\(G(u)=\frac{\sqrt{(1-(u-1)t)^2-4t} } {1-(u-1)^2 t} \) is a power
series in both \(u\) and \(t\), we have
\begin{align*}
\CT_x \frac{1-x}{1+x}\frac{\Delta}{D_2/x^2} &=\CT_u
\frac{\sqrt{(1-(u-1)t)^2-4t} } {1-(u-1)^2
t}=\frac{\sqrt{(1+t)^2-4t}}{1-t}
= 1.
\end{align*}

c) By \eqref{P1} and \eqref{P2}, the following easily verified
equation (from later calculation)
\begin{multline*}
\frac{\partial}{\partial
t}\frac{(1-x)(1-3tx)}{(1-t)^2}\Delta=\frac{1+3t}{(t-1)^3}(x-x^{-1})\Delta\\
-x\frac{\partial}{\partial
x}\frac{4tx^3+(3t^2-7t)x^2+(3t^2-4t-3)x+3t^2+t}{2(t-1)^3tx}\Delta
\end{multline*}
shows that \(\CT_x(1-x)(1-3xt)(1-t)^{-2}\Delta\) is a constant.
Equation \eqref{P8} thus follows by checking the \(t=0\) case.
\qed

\begin{rem}
In Lemma \ref{eliminat-lem}, part a) can also be regarded as
applications of Jacobi's residue formula by letting \(y=x^{-1}\).
We suspect that Jacobi's residue formula can also be used to prove
part c), which arises naturally when proving the differential
equation for \(\mathcal{C}(t)\). See \cite[Proposition
1]{Xin-3noncrossing}.
\end{rem}

{\bf{Proof of Proposition \ref{laurent-prop}.}}  
We will successively reduce \(p,q,r,s\) to $0$, so it is sufficient
to deal with the  cases of $p,q,r,s\geq 1$.

Let \(\deg_x D\) be the degree of \(D\) in \(x\). By classical
results for partial fraction decompositions, we have the unique
decomposition
\[\frac{L(x,t)}{(1+x)^{p}D_1^qD_2^rD_3^s} = l(x,t)+\sum_{i=1}^p
\frac{P_i }{D_0^i}+\sum_{i=1}^q \frac{Q_i(x,t)
}{D_1^i}+\sum_{i=1}^r \frac{R_i(x,t) }{D_2^i}+\sum_{i=1}^s
\frac{S_i(x,t) }{D_3^i},
\]where \(l(x,t)\) is a Laurent polynomial, \(D_0=1+x\), \(P_i\in
\mathbb{R}\), \(\deg_x Q_i(x,t)<\deg_x D_1\), \(\deg_x
R_i(x,t)<\deg_x D_2\) and \(\deg_x S_i(x,t)<\deg_x D_3\) for all
\(i\).

 \vspace {1mm}
We shall often use the the above decomposition when multiplied
through by \(\Delta\), so we are actually dealing with a
$\mathbb{Q}(t)$-linear combination of $x^k\Delta/D_i^j$, where $0\le
k \le \deg_x D_i$ if $j\ge 1$ and $k\in \mathbb{Z}$ if otherwise.
Let us call \(x^k\Delta/D_i^{j}\) together with its coefficient the
\emph{\(x^k\Delta/D_i^{j}\)-term}, and the collection of
\(x^k\Delta/D_i^{j}\)-terms for $0\le k\le \deg_x D_i$ the
\emph{\(x^*\Delta/D_i^{j}\)-term}. We will subtract by known
constant terms to reduce our original constant term to simpler
forms.

\vspace {1mm} Step 1:  Reduce $p,q,r$ to $1$ by the following
procedure. Successively eliminate the \(x^*\Delta/D_2^{r}\)-term,
and then the \(x^*\Delta/D_2^{r-1}\), \dots,
\(x^*\Delta/D_2^{2}\)-terms, and similarly for \(D_1\) and \(D_0\).
The process works for any irreducible polynomial \(D=D(x,t)\) that
is coprime to \(x\) and \(D_3\). {Denote by \(N\Delta/D^r\) the
\(x^*\Delta/D^{r}\)-term where \(N=N(x,t)\) and \(\deg_xN <
\deg_xD\).}
  Noticing
\[D_3=x^2\Delta^2 \Rightarrow \frac{\partial \Delta}{\partial
x}=-\frac{\Delta}{x}+\frac{1}{2}\frac{\Delta}{D_3}\frac{\partial
D_3}{\partial x},\] we can eliminate   \(N\Delta/D^r\)  for \(r\ge
2\)   by subtracting the partial fraction decomposition of the
following constant term.
\begin{multline}0=\CT_x x\frac{\partial }{\partial x} \frac{S
\Delta}{D^{r-1}}
=\CT_x\left( \frac{x
 \Delta}{D^{r-1}}\frac{\partial S}{\partial x}-\frac{S\Delta}{D^{r-1}}+\frac{1}{2}
 \frac{xS\Delta}{D^{r-1}D_3}\frac{\partial D_3}{\partial x}
 \right.\\
  \left.
 +\frac{(1-r)xS\Delta}{D^r}\frac{\partial D}{\partial x}\right).\label{forD}\end{multline} Here \(S\) is
an appropriately chosen polynomial in \(x\) such that \(D\) divides
\((1-r)\frac{\partial D}{\partial x}xS-N \). Since \(D\) is
irreducible and coprime to \(x\), it is coprime to \(x\frac{\partial
D}{\partial x}\). Therefore we can find polynomials \(\alpha\) and
\(\beta\) in \(x\) (by the Euclidean algorithm) such that
\begin{align}\alpha D+\beta x\frac{\partial D}{\partial
x}=1.\label{prime-eq}\end{align} Now choose
\[S=N\beta/(1-r)\Rightarrow (1-r)\frac{\partial D}{\partial
x}xS-N=-\alpha ND.\]

 Step 2: Reduce $p$ and $q$ to $0$. First eliminate the \(x^*\Delta/D_0\)-term by using \eqref{P5},
which can be rewritten as
\[\CT_x(1+x)^{-1}\Delta=(1-t)/2+1/2\CT_x\Delta.\]
Next eliminate the \(x^*\Delta/D_1\)-term  by subtracting a linear
combination of the following two constant terms.
\begin{equation*}
  \CT_x
  \frac{1-x}{1+x}\frac{x\Delta}{D_1}-1=\CT_x\left(-\frac{\Delta}{t-1}-\frac{4t+(5t-1)x}{(t-1)D_1}\Delta\right)=0,
  \end{equation*}
  \vspace{-6mm}
\begin{multline*}\CT_xx\frac{\partial }{\partial x}\Big(\ln(1-Y/x)-\frac{1}{2}\ln(D_1/x)\Big)\\
=\CT_x\Big(-\frac{1}{4}+\frac{\Delta}{4(t-1)}+\frac{x\Delta}{D_1}
- \frac{t^2+(t^2-2t)x+(t^2-t-1)x^2+tx^3}{2(t-1)D_3}\Delta\Big)=0.
\end{multline*}
Step 3: Eliminate all the \(x^k\Delta/D_2\)-terms for \(k=1,2,3\) by
using the following three constant terms.
\begin{equation*}
 \CT_x  \frac{(x^3-x){\Delta}}{D_2}=0 \ \ \ \ \ \ \ \   (\textrm{by \eqref{P4} with } k=1)
 \end{equation*}
 \vspace{-6mm}
\begin{multline*}
 \CT_x\frac{1-x}{1+x}\frac{x^2\Delta}{D_2}-1=\CT_x\left(-\frac{\Delta}{t-1}-\frac{2t+4tx+(3t-1)x^2+2tx^3}{(t-1)D_2}\Delta\right)=0
 \end{multline*}
  \vspace{-6mm}
 \begin{multline*}
 \CT_xx\frac{\partial }{\partial
x}\Big(\ln(1-xY)-\frac{1}{2}\ln(D_2/x^2)\Big)\\
=\CT_x\left(\frac{3}{4}+\frac{\Delta}{4(t-1)} -\frac{\Delta
x^2}{D_2}
 -\frac{t^2+(t^2-2t)x+(3t^2-5t+1)x^2+tx^3}{2(t-1)D_3}\Delta\right)=0
\end{multline*}

Step 4: Reduce the current $s$ to $0$. Eliminate one by one (if
needed) the \(x^*\Delta/D_3^{\ell}\)-terms for $\ell=s,s-1,\dots, 1$
similarly as in Step 1.
By collecting terms in \eqref{forD} (with \(D=D_3\)), 
with \(\alpha\) and \(\beta\) in \eqref{prime-eq}, we can eliminate
 \(N\Delta/D_3^r\)  by choosing \[S=N\beta/(3/2-r) \Rightarrow
\big(\frac{3}{2}-r\big)\frac{\partial D_3}{\partial x}xS-N=-\alpha
ND_3.
\]

Step 5: Remove all of the \(x^k\Delta\)-terms for $k\le 0$ or $k\ge
2$. First eliminate all \(x^k\Delta\)-terms for \(k<0\) by
\eqref{P2}. Then eliminate all \(x^k\Delta\)-terms for
\(k=\ell,\ell-1, \dots,3\), where
\(\ell=\max\{\deg_xl(x,t),\deg_xl(\xx,t)\}\), one
by one by the formulas 
\begin{align*}
x\frac{\partial}{\partial x}x^{1+i} \Delta^{3}&= \big((4+i)t^2
x^{i+3}+b
x^{i+2}+b'x^{i+1}+b''x^i\big)\cdot\Delta, \text{ for } i\geq 1,\\
x\frac{\partial}{\partial
x}x\Delta^{3}&=\left(4t^2x^3+(5t^2-5t)x^2+(3t^2+1-6t)x-t^2+t-2\frac{t^2}{x}\right)\cdot
\Delta,
\end{align*}
in which the \(b\)'s are independent of \(x\). Finally eliminate
\(x^k\Delta\)-terms  by \eqref{P8} for \(k=2\), and by \eqref{P2}
again for \(k=-1\).

Step 6: Reduce $r$ to $0$. By Step 3, it is sufficient to eliminate
the \(x^0\Delta/D_2\)-term. This is done by showing the following
equality:
\begin{multline}\CT_x \Big({\frac {9}{32}} \,{\frac{1}{{t}^{2} ( 9\,t-1 )
}}(-1+8\,t+55\,{t}^{2}-440\,{t}^{3}+861\,{t}^{4}-528\,{t}^{5}+45\,{
t}^{6})\\
+\big({\frac {9}{32}}\,{\frac{1}{{t}^{2} ( 9\,t-1 )
}}(1-5\,t-74\,{t}^{2}+210\,{t}^{3}-87\,{t}^{4}-45\,{t}^{5})\\
+{\frac {9}{32}}\,{\frac{1}{{t}^{2} ( 9\,t-1 )
}}(-4+8t+240t^2-552t^3+180t^4)tx\big)\Delta\\
-{\frac {9}{8}}\, { \frac {\Delta\, (1+5t-21t^2+15t^3 ) }{ t D_2
}}\Big)=0.\label{forD2}\end{multline}

Denote by \(E(t)\) the the left-hand side of the above equation. To
show that  \(E(t)=0\), we first show that \(E(t)\) satisfies a
\(D\)-finite equation. The method is typical.

Using Steps 1--5, we can rewrite \[\frac{d^i}{d
t^i}E(t)=\widetilde{R}_i(t)+\CT_x(\widetilde{P}_i(t)+\widetilde{Q}_i(t)x)\Delta+\CT_x\widetilde{S}_i(t)\Delta/D_2,
\ \ \ i=0,1,2,3,4.\] By solving the system of equations
\begin{equation*}
\left\{
\begin{array}{c}
a\widetilde{R}_0(t)+b\widetilde{R}_1(t)+c\widetilde{R}_2(t)+d\widetilde{R}_3(t)+e\widetilde{R}_4(t)=0\\
a\widetilde{P}_0(t)+b\widetilde{P}_1(t)+c\widetilde{P}_2(t)+d\widetilde{P}_3(t)+e\widetilde{P}_4(t)=0\ \ \\
a\widetilde{Q}_0(t)+b\widetilde{Q}_1(t)+c\widetilde{Q}_2(t)+d\widetilde{Q}_3(t)+e\widetilde{Q}_4(t)=0\\
a\widetilde{S}_0(t)+b\widetilde{S}_1(t)+c\widetilde{S}_2(t)+d\widetilde{S}_3(t)+e\widetilde{S}_4(t)=0\
\
\end{array}\right.
\end{equation*}
for \(a,b,c,d,e\) independent of \(x\), we get the nontrivial
solution
\[a=1,b=\frac{2t+10t^2-42t^3+30t^4}{3+5t+21t^2-45t^3},c=d=e=0.\]
This implies that
 \(E(t)+b\cdot
\frac{d }{d t}E(t)=0\).   Solving this differential equation gives
\[E(t)=C_0(15t^{\frac{3}{2}}-21t^{\frac{1}{2}}+5t^{-\frac{1}{2}}+t^{-\frac{3}{2}})\]
for some constant \(C_0\).

On the other hand, by using {\sc Maple} to expand \(E(t)\) as a
series in \(t\) and then take constant term in \(x\), we see that
\(E(t)\) is actually a power series in \(t\) with \(E(0)=0\). It
then follows that \(C_0\) must be \(0\) and hence \(E(t)=0\) as
desired.   \qed

Once Proposition \ref{laurent-prop} is established, the differential
equations (e.g., \eqref{deqfornco}) can be proved by {\sc Maple}.
The package can be downloaded at \\
 \emph{
http://www.combinatorics.net.cn/homepage/xin/maple/bs3np.txt}.

\section{Analogous Results for Bi-symmetric Enhanced 3-noncrossing
Partitions}\label{sec:hestating}

Chen et al. \cite{chen} also considered a variation of $k$-crossings
(nestings), called enhanced $k$-crossings (nestings). Given a
partition $P$ of $[n]$, its enhanced graph representation is
obtained by adding a loop to each isolated point in the graph
representation of $P$. Then an enhanced $k$-crossing of $P$ is a set
of $k$ edges $(i_1, j_1), (i_2, j_2),\ldots , (i_k, j_k)$ of the
enhanced representation of $P$ such that $i_1<i_2< \cdots< i_k \leq
j_1<j_2<\cdots<j_k$. Our approach for counting bi-symmetric
3-noncrossing partitions can be easily adapted to obtain analogous
enumeration results for bi-symmetric partitions avoiding enhanced
3-crossings.

\medskip
Let \(\widetilde{E}_3(n)\) be the number of bi-symmetric partitions
of \([n]\) avoiding enhanced 3-crossings. We obtain the following
result.
\begin{prop}\label{bisym-enhanced}
The numbers \(\widetilde{E}_3(2n)\) satisfy
\(\widetilde{E}_3(0)=1\), \(\widetilde{E}_3(2)=2\), and
\begin{multline*}
8(n+3)(n+1)\widetilde{E}_3(2n)+(7n^2+41n+58)\widetilde{E}_3(2n+2)-(n+4)(n+5)\widetilde{E}_3(2n+4)=0.
\end{multline*}
The numbers \(\widetilde{E}_3(2n+1)\) satisfy
\(\widetilde{E}_3(1)=1\), \(\widetilde{E}_3(3)=3\),
\(\widetilde{E}_3(5)=11\),   and
\begin{multline*}
32(n+2)^2\widetilde{E}_3(2n+1)+(36n^2+220n+328)\widetilde{E}_3(2n+3)\qquad
\qquad
\\
+(3n^2+26n+56)\widetilde{E}_3(2n+5)-(n+6)^2\widetilde{E}_3(2n+7)=0\end{multline*}
Equivalently, their associated generating functions
\(\H_e(t)=\sum_{n\geq0}\widetilde{E}_3(2n)t^n\) and
\(\H_o(t)=\sum_{n\geq0}\widetilde{E}_3(2n+1)t^n\) satisfy
\vspace{-5mm}
\begin{multline*}6-(6-24t-24t^2)\H_e(t)-(6t-34t^2-40t^3){\frac {d}{d{t}}}\H_e(t)
\\
+t^2(t+1)(8t-1){\frac {d^{2}}{d{t}^{2}}}\H_e(t)=0,\end{multline*}
\vspace{-10mm}
\begin{multline*}
(9+32t+32t^2)+(-9+16t+144t^2+128t^3)\H_o(t) +
                  t(-7+17t\\ +184t^2+160t^3){\frac
{d}{d{t}}}\H_o(t) +(4t+1)(8t-1)(t+1)t^2{\frac
{d^{2}}{d{t}^{2}}}\H_o(t)=0.
\end{multline*}
\end{prop}

We need the lattice walk interpretations. A hesitating lattice walk
satisfies the following walking rules: when pairing every two steps
from the beginning, each pair of steps has one of the following
three types: i) a stay step followed by an $e_i$ step, ii) a  $-e_i$
step followed by a stay step, iii) an $e_i$ step followed by a
$-e_j$ step.  It was  pointed out that partitions of \([n]\)
avoiding enhanced \(k+1\)-crossings are in bijection with hesitating
tableaux of height bounded by \(k\) under a map \(\bar \phi\) in
\cite{chen}. In turn, these hesitating tableaux are in one-to-one
correspondence with certain $W_k$-hesitating lattice walks. For the
$k=2$ case, this reduces  to a bijection  between  partitions of
\([n]\) avoiding enhanced 3-crossings and $W_2$-hesitating lattice
walks of length $2n$ starting and ending at the point $(1,0)$.

Given a set $A$ of points,   let $E^A(n)$ be the number of
$W_2$-hesitating lattice walks of length $2n$ starting from points
in $A$ and ending at $(1,0)$, and let $\mathcal{E}^A(t)$ be the
generating function of $E^A(n)$. Similar approach as for the
vacillating case can give us the following analogous result.
\begin{prop}\label{basic-prop-en}
For any set $A$ of lattice points in \(W_2\), we have
\begin{equation*}
\mathcal{E}^A(t)=\CT_x\frac{(x^{-2}-x^3)\big(x\widetilde{Y}A(x,\widetilde{Y})-x^{-1}\widetilde{Y}^2A(x^{-1}\widetilde{Y},\widetilde{Y})+x^{-2}\widetilde{Y}A(x^{-1}\widetilde{Y},x^{-1})\big)}{t(1+x)},
\end{equation*}
where $\widetilde{Y}=\widetilde{Y}(x;t)$ is the unique power
series in $t$ satisfying
$\widetilde{Y}=t(1+x^{-1})(1+\widetilde{Y})(x+\widetilde{Y})$
given by
$$\widetilde{Y}=\frac{1-tx^{-1}(1+x)^2-\sqrt{((1-tx^{-1}(1+x)^2)^2-4t^2x^{-1}(1+x)^2}}{2(1+x^{-1})t}.$$
\end{prop}

Again, by the    correspondence \(\bar \phi\),   a partition $P$ of
$[n]$ is bi-symmetric if and only if the corresponding hesitating
lattice walk is palindromic. By a parallel argument as for the
vacillating case, and observing that the $n+1$st pair of steps for
each palindromic hesitating lattice walk of length $4n+2$ must be an
$e_i$ step followed by a $-e_i$ step for some $i$, we can obtain
formulas for $\H_e(t)$ and $\H_o(t)$:
\begin{align*} \H_e(t) &=\E^{A_2(x,y)}(t),\\
\H_o(t)&=\E^{A_4(x,y)+2(A_2(x,y)-A_4(x,y))}(t)=\E^{2A_2(x,y)-A_4(x,y)}(t).
\end{align*}

The above observations and Proposition \ref{basic-prop-en} enable us
to develop a similar {\sc{Maple}} package for 2-dimensional
hesitating lattice walks enumerating problems. Actually we can use
our package for the vacillating case by redefining some initial
variables. See Appendix \ref{appendix-hesitating}. With our package,
we can prove Proposition \ref{bisym-enhanced} in a second. Moreover,
we find the following result.
\begin{prop}\label{en-Baxter}The numbers $E^{A_3}(n)$ satisfy $E^{A_3}(0)=1$,
$E^{A_3}(1)=2$, and
\begin{multline*}
8(n+1)(n+2)E^{A_3}(n)+(7n^2+49n+82)E^{A_3}(n+1)\\
-(n+5)(n+6)E^{A_3}(n+2)=0.
\end{multline*}
Equivalently, its associated generating function satisfies that
\begin{multline*}12+4(-3+10t+4t^2)\mathcal{E}^{A_3}(t)+2t(-4+21t+16t^2)\frac{d}{dt}\mathcal{E}^{A_3}(t)\\
+t^2(t+1)(8t-1)\frac{d^2}{dt^2}\mathcal{E}^{A_3}(t)=0.\end{multline*}
\end{prop}
By searching through \cite[A001181]{Sloane}, we discover that the
number of hesitating lattice walks of length $2n$ starting from
$(1,0)$ and ending in $A_3$ is equal to the number $b_{n+1}$ of
Baxter permutations of length $n+1$. To prove it, we use the
formula
$$b_n=\frac{2}{n(n+1)^2}\sum_{k=0}^{n-1}{n+1 \choose k}{n+1 \choose k+1}{n+1 \choose k+2},$$
and apply the creative telescoping of \cite{AB}. It is worth
mentioning that $b_n$ also counts the number of watermelons
consisting of three vicious walkers. See \cite{Melou-oscu} and
\cite{Dulucq}. Note that there are 8 possible pair of steps for
$W_2$-hesitating lattice walks, and 8 possible 1-steps for
watermelons consisting of three vicious walkers. Then a natural
question arises: Can we find a bijection between them?

Let $\tilde{w}_2((1,0),A,n)$ be the number of $W_2$-hesitating
lattice walks of length $n$, starting at $(1,0)$ and ending in
$A$. We conclude this subsection by Table \ref{tab3} of some
asymptotic estimates.

\medskip

\begin{figure}[ht]
\begin{center}
\begin{tabular}{c|llllllllllll}
n&0&1&2&3&4&5&6&\(\rightarrow\)&\(\infty\)\\
\hline
\(\tilde{w}_2((1,0),A_1,2n)\)&1& 1& 2& 5& 15& 51& 191& \(\sim\)& \(\lambda_1\cdot8^n/n^7\)\\
\(\tilde{w}_2((1,0),A_2,2n)\)&1&2&7&29&136&692&3739&\(\sim\)&\(\lambda_2\cdot 8^n/n^3\)\\
\(\tilde{w}_2((1,0),A_2,2n+1)\)&1&3&11&48&232&1207&6631&\(\sim\)&\(\lambda_3 \cdot 8^n/n^3\)\\
\(\tilde{w}_2((1,0),A_3,2n)\)&1&2&6&22&92&422&2074&\(\sim\)&\(\lambda_4\cdot
8^n/n^4\)

\end{tabular}

\caption{  The first several numbers of hesitating lattice walks and
their asymptotic estimate, where $\lambda_1\approx
6670.312,\lambda_2\approx 7.835,\lambda_3\approx 15.669$, and
$\lambda_4\approx 46.988$. \label{tab3} }
\end{center}
\end{figure}

\section{Discussion}
Since our discussion for the vacillating case and that for    the
hesitating case are similar to each other, we focus on the
vacillating case.

The very general theory in \cite{lipshitz-df} asserts that
\(\mathcal{C}^{A}(t)\) is \emph{D-finite} if \(A(x,y)\) is rational.
That is, it satisfies a linear differential equation with polynomial
coefficients, or equivalently, \(C^A(n)\) satisfies a P-recurrence.
However the degree of the equations suggested in \cite{lipshitz-df}
is usually too large for proving simple P-recurrences as we
consider. Note that these recurrences can be easily guessed, using
the {\sc Maple} package {\sc Gfun}. The recurrence for \(C_3(n)\)
was proved by using the Lagrange inversion formula to give a single
sum formula and then applying the creative telescoping of \cite{AB}.
However, the same route is difficult to apply to our case. The
Lagrange inversion formula will give us a complicated double sum.

Actually our {\sc Maple} package can produce the differential
equation for \(\C^A(t)\) for any \(A=P(x,y)/((1-x)(1-xy))\) with
\(P(x,y)\) a polynomial. The whole process will be completed within
seconds if \(P(x,y)\) is simple. Two curious observations are worth
mentioning. We have described how to write \(\C^A(t)\) as \(\CT_x
T_0+\CT_x T_1\Delta\) for rational \(T_0\) and \(T_1\), and
Proposition \ref{laurent-prop} deals with \(\CT_x T_1\Delta\). In
practice, we find that i) \(T_0\) does not contain \(D_2\) and
\(D_3\) as denominators; ii) using the constant term identity
\[
\CT_x\frac{4t+(5t-1)x}{2t+(5t-1)x+2tx^2}=1,  \]  obtained by
considering the constant term of \(x\frac{\partial}{\partial x}\ln
D_1\), one sees that  \(\CT_x T_0\) is always   a rational function
in \(t\). We do not know why \(\CT_x T_0\) is always rational, since
this is not true if we take, e.g., \(A(x,y)=x/(1-x)^2\) or if we
pick out a term from the sum in \eqref{ctfunction}. It is not a
problem even if \(T_0\) has \(D_1,D_2,D_3\) as denominators. We can
suitably enlarge \(\mathcal{L}\) and increase the dimension bound to
fit in our package.

In the proof of Proposition \ref{laurent-prop}, only using Steps
1,4,5, one can already give an upper bound for the dimension of
\(\mathcal{L}\). Theoretically one can prove differential equations
like \eqref{deqfornce} similarly as in Step 6. Such equations, once
proved, will reduce the upper bound of the dimension. Equation
\eqref{P8} is actually obtained when proving the differential
equation satisfied by the generating function $\C(t)$ (See
\cite[Proposition 1]{Xin-3noncrossing}); The equation \(E(t)=0\) is
obtained when proving \eqref{deqfordiag}. Finding small upper bounds
for this type of problems may help discovering and proving new
formulas, and possibly reducing the upper bound again. This idea may
well apply to other situations.

Our contribution is to reduce the upper bound to only 3. This
results in a fast algorithm for  2-dimensional vacillating lattice
walk  enumeration problems. The number 3 should be the actual
dimension of \(\mathcal{L}\), since otherwise \(\C(t)\) must satisfy
a lower degree differential equation, which is not suggested by the
{\sc Maple} package {\sc Gfun}. However it seems hard to prove the
equality.

 Moreover, it
would be interesting to find some    combinatorial proofs  for the
interesting  relations  stated in Propositions \ref{twice-relation},
  \ref{vaci-odd} and   \ref{en-Baxter}.

{\bf{Acknowledgments}}

The authors would like to thank the referees for helpful
suggestions to improve the presentation, and Christian
Krattenthaler and Tom Roby for valuable comments. This work was
supported by the 973 Project, the PCSIRT project of the Ministry
of Education, the Ministry of Science and Technology and the NSF
of China.

\appendix

\section{Appendix: Enumeration of Bi-symmetric Noncrossing Partitions
\label{sec:noncrossing}} In this section we consider the enumeration
of bi-symmetric noncrossing partitions. To state our result, we need
the following definition: An \emph{oscillating tableau} (or
\emph{up-down tableau}) of shape \(\mu\) and length \(n\) is a
sequence \((\emptyset=\mu^0,\mu^1,\dots,\mu^{n}=\mu)\) of partitions
such that for all \(1\le i\le n -1\), the diagram of \(\mu^i \) is
obtained from \(\mu^{i-1}\) by either adding or removing one square.
\begin{prop}\label{propfornp}
There is a bijection between the set of palindromic   oscillating
tableaux of length \(2n\) and  height bounded by \(1\) and
 the set of palindromic  vacillating
tableaux of length \(2n\) and  height bounded by \(1\). Moreover,
both of them are enumerated by \({n\choose \lfloor n/2 \rfloor}\).
\end{prop}

To construct the bijection, it is convenient to introduce an
intermediate set \(\mathcal{W}(n)\) of all \(01\) words of length
\(n\) with no initial segments containing more 0's than 1's. It is
well-known that \(|\mathcal{W}(n)|={n\choose \lfloor n/2 \rfloor}\).
For any \(01\) word \(w\) of length \(n\) and \(s=\)0 or 1, define
\vspace{-3mm}
\[odd(w,s)=|\{i\ \mathrm{is\ odd}\mid w_i=s, \ i\in[n] \}|,\]
\vspace{-8mm}
\[even(w,s)=|\{i\ \mathrm{is\ even}\mid w_i=s, \ i\in[n] \}| .\] Then
we have the following characterization.
\begin{lem}\label{lem-word}
 \(w\in \mathcal{W}(n)\Leftrightarrow \mathrm{for\ any\ initial\ segment\ }w'\mathrm{\ of\ }w\mathrm{,\ }  {even}(w',1)\geq {odd}(\\w',
 0)\).
\end{lem}
{\bf{Proof.}} Let \(w'\) be the initial segment of \(w\) of length
\(m\). Then we have the natural equality
\begin{equation}
even(w',1)+even(w',0)+\chi(m  \mathrm{\ is\ odd})=
odd(w',0)+odd(w',1),\label{e-word1}
\end{equation}
where \(\chi(S)\) is \(1\) if the statement \(S\) is true and \(0\)
otherwise. On the other hand, by definition \(w\in \mathcal{W}(n)\)
if and only if for every initial segment \(w'\) of \(w\) we have
\begin{equation}
even(w',1)+odd(w',1)\geq even(w',0)+odd(w',0).\label{e-word2}
\end{equation}
Obviously \eqref{e-word2} can be replaced with
\(\eqref{e-word1}+\eqref{e-word2}\), which is equivalent to
\({even}(w',1) \geq {odd}(w',0) \). \qed

{\bf{Proof  of Proposition \ref{propfornp}.}}  Given a palindromic
oscillating tableau \(O=(O_0,O_1, \dots,O_{2n})\) of length \(2n\)
and  height  bounded by \(1\), we have a natural encoding
\(\theta(O)=w=w_1w_2\cdots w_n \in \mathcal{W}(n)\) defined by
\(w_i=1\) if \(O_i\) is obtained from \(O_{i-1}\) by adding a
square, and \(w_i=0\) otherwise. Note that palindromic means that
\((O_0,O_1,\dots,O_n)\) already carries all information of \(O\).

Next we conclude  the proposition by constructing a bijection
\(\eta\) from the set of palindromic vacillating tableaux of length
\(n\) and  height  bounded by \(1\) to \(\mathcal{W}(n)\). Given
such a tableau \(V=(V_0,V_1,\ldots,V_{2n})\), we define
\(\eta(V)=w=w_1w_2\cdots w_n\) according to the four cases:  (i) if
\(i\) is odd and \(V_i=V_{i-1}\), then \(w_i=1\); (ii) if \(i\) is
even and \(V_i\) is obtained from \(V_{i-1}\) by adding a square,
then \(w_i=1\); (iii) if \(i\) is even and \(V_i=V_{i-1}\), then
\(w_i=0\); (iv) if \(i\) is odd and \(V_i\) is obtained from
\(V_{i-1}\) by deleting a square, then \(w_i=0\). Clearly \(V\) is a
vacillating tableaux if and only if the number of type (ii) moves is
no less than the number of type (iv) moves in any initial segment of
\(V\). This is the same as that in any initial segment \(w'\) of
\(w\), \(even(w',1)\geq odd(w',0)\), which is equivalent to \(w\in
\mathcal{W}(n)\) by Lemma \ref{lem-word}. Thus \(\eta\) is the
desired bijection. \qed

\medskip
\noindent {\bf Example.} Let \(O=(0,1,2,1,2,1,0,\dot{1},
0,1,2,1,2,1,0)\) be the palindromic oscillating tableau, where the
integers stand for one row partitions and we put a \(\cdot\) over
the cental diagram. Then \(\theta(O)=1101001\), and the
corresponding palindromic vacillating tableau is
\((0,0,1,0,1,1,2,\dot{1},2,1,1,0,1,0,0)\).

\begin{rem}
A word \(w\in \mathcal{W}(2n)\) consisting of \(n\) 1's and \(n\)
0's is called a {\it{Dyck word}}. Denote the set of such words by
\(\mathcal{D}(n)\). By the proof of Lemma \ref{lem-word}, we observe
that \[w\in \mathcal{D}(n) \Leftrightarrow even(w,1)=odd(w,0), \quad
w \in \mathcal{W}(2n).\] When the   bijections \(\theta\) and
\(\eta\) are restricted to \(\mathcal{D}(n)\), we can obtain a
bijection between noncrossing matchings of \([2n]\) and noncrossing
partitions of \([n]\).
\end{rem}

\section{Appendix: Initial Variables for Hesitating Lattice Walks \label{appendix-hesitating}}
To apply the vacillating case package to the hesitating case, we
reset the initial variables as follows:
\[u=(\xx+2+x)=\xx(1+x)^2, \]
 \[\Delta\equiv\Delta(x,t)=\sqrt {((1-tx^{-1}(1+x)^2)^2-4t^2x^{-1}(1+x)^2}=\sqrt {(1-ut)^2-4ut^2},\]
\[D_1= x-2t(1+x)^2=x(1-2tu),\]
\[D_2= -t-2tx-(2t-1)x^2-2tx^3-tx^4=x^2(1+t-t(1-u)^2),\]
\[D_3=t^2-2tx-(2t^2+4t-1)x^2-2tx^3+t^2x^4=x^2((1-ut)^2-4ut^2).\]

Similarly,   $D_3=x^2\Delta^2$ and
\begin{align*}
 \frac{\partial }{\partial x}
\Delta(x,t)&= \frac{t^2x^4-tx^3+tx-t^2}{xD_3}\Delta,
 \\
\frac{\partial }{\partial t} \Delta(x,t)&
=\frac{tx^4-x^3-2tx^2-2x^2-x+t}{D_3}\Delta.
\end{align*}

The following is a replacement of Lemma \ref{eliminat-lem}.
\begin{lem}
a) For all \(k\in \mathbb{Z}\), we have \vspace {-3mm}
\begin{align*}
\CT_x(x^k-x^{-k})\Delta&=0,\\
\CT_x   (x^k-x^{-k})\frac{x\Delta} {D_1}&=0,\\
\CT_x   (x^k-x^{-k})\frac{x^2\Delta}{D_2}&=0.
\end{align*}
b)\vspace {-8mm} \begin{align*}
\CT_x\frac{1-x}{1+x}\Delta&=1.\\
\CT_x \frac{1-x}{1+x}  \frac{x\Delta}{D_1}&=1.\\
\CT_x\frac{1-x}{1+x} \frac{x^2\Delta}{D_2}&=1.
 \end{align*}
c)\vspace {-6mm} \begin{align*}\CT_x
(1-x)(1+t-3xt)\Delta=1.\end{align*}
\end{lem}
The proofs of part a) and b) are similar. Part c) follows from the
following equality
\begin{multline*}
\frac{\partial}{\partial
t}(1-x)(1-3tx+t)\Delta=-\frac{1+4t}{1+t}(x-x^{-1})\Delta\\
-x\frac{\partial}{\partial
x}\frac{-4t(t+1)x^3+t(4t+7)x^2+(4t^2+10t+3)x-(1+4t)t}{2t(1+t)x}\Delta.
\end{multline*}

The following equality is an analogy of \eqref{forD2}.
\begin{multline*}\CT_x
\Big(\frac{3}{32}\frac{1-6t-32t^2+48t^3+64t^4}{t^2(t+1)^2(8t-1)} +\big(\frac{3}{32}\frac{1-4t-44t^2-48t^3}{t^2(1+t)^2(1-8t)}\\
-\frac{3}{8}\frac{1-2t-28t^2-16t^3}{t(1+t)^2(1-8t)}x\big)\Delta +\,
{
\frac{3}{8}\frac{(1+4t)\Delta}{t(1+t)D_2}}\Big)=0.\label{forD2}\end{multline*}



\begin{thebibliography}{00}




\bibitem{bousquet-petkov}
M.~Bousquet-M{\'e}lou and  M. Petkov\v sek,
  Walks confined in a quadrant are not always $D$-finite,
 {Theoret. Comput. Sci.}  \textbf{307} (2003), 257--276.



\bibitem{bousquet-motifs}
M.~Bousquet-M\'elou,
  Four classes of pattern-avoiding permutations under one roof:
  generating trees with two labels,
 {Electron.~J.~Combin.}  \textbf{9} (2003), R19.

\bibitem{bousquet-kreweras}
M.~Bousquet-M{\'e}lou,
  Walks in the quarter plane: Kreweras' algebraic model,
 { Ann. Appl. Probab.}  \textbf{15} (2005), 1451--1491.


\bibitem{Xin-3noncrossing}
M. Bousquet-M\'elou and G. Xin, On partitions avoiding
3-crossings,  {S\'em. Lothar. Combin.}  {\bf 54} (2005), Art.
B54e.

\bibitem{Melou-oscu}  M. Bousquet-M\'{e}lou, Three osculating
walkers,  {J. Phys.: Conf. Ser.}  \textbf{42} (2006), 35--46.

\bibitem{chen}W. Y. C. Chen, E. Y. P. Deng, R. R. X. Du, R. Stanley and  C. H. Yan, Crossings and
nestings of matchings and partitions,     {Trans. Amer. Math.
Soc.}  \textbf{359} (2007), 1555--1575.

\bibitem{Graham-Chung} F. R. K. Chung, R. L. Graham, V. E. Hoggatt Jr and M. Kleiman,
  The number of Baxter permutations, {J. Combin. Theory Ser.
A}  \textbf{24} (1978),  382--394.

\bibitem{Dulucq} S. Dulucq  and  O.  Guibert,   Baxter permutations,    {Discrete Math.}  \textbf{180}
(1998), 143--156.

 \bibitem{Fomin-GRSKC} S. V. Fomin,  Generalized Robinson-Schensted-Knuth
 correspondence,   {J. Soviet Math.}   {\bf 41}   (1988),
 979--991.

  \bibitem{Fomin-Sadgg} S. V. Fomin,  Schensted algorithms for dual graded graphs,   {J. Algebraic Combin.}  {\bf 4}   (1995),
   5--45.


 \bibitem{Fomin-Sk} S. V. Fomin,  Schur operators and Knuth correspondences,  {J. Combin. Theory Ser. A}   {\bf 72} (1995),
 277--292.






 \bibitem{Knuth}
 D. E. Knuth, Permutations, matrices and generalized Young tableaux,  {Pacic J.
 Math.}
\textbf{34} (1970), 709--727.


\bibitem{Krattenthaler} C. Krattenthaler, Growth diagrams, and increasing and decreasing
chains in fillings of Ferrers shapes,  {Adv. in Appl. Math.}
\textbf{37} (2006),
 404--431.

\bibitem{lipshitz-df}
L.~Lipshitz,
  D-finite power series,
 { J. Algebra}   \textbf{122} (1989), 353--373.

\bibitem{AB}
M.~Petkov{\v{s}}ek, H.~S. Wilf  and D.~Zeilberger,
\newblock {\em {\(A=B\)}},
\newblock A K Peters Ltd., Wellesley, MA, 1996.

%
%
\bibitem{Roby} T. W. Roby, Applications and Extensions of Fomin's Generalization of the Robinson-Schensted Correspondence
to Differential Posets, PhD thesis, MIT, Cambridge, MA, 1991.

\bibitem{SLR} M. P. Sch\"utzenberger,    La corres\-pondance de Robinson,
  in ``Combinatoire et Repr\'esentation du
  Groupe Sym\'etrique",  pp. 59--113. Lecture Notes in Math., Vol. 579, Springer, Berlin,
  1977.




\bibitem{Sloane} N. J. A. Sloane, The On-Line Encyclopedia of Integer
Sequences, published electronically at
 \emph{www.research.att.com/\raisebox{-3pt}{\~{}}njas/sequences}.


\bibitem{Rstanley} R. P. Stanley, Enumerative Combinatorics, vol. 2,
Cambridge University Press, Cambridge, 1999.


%






\bibitem{Wimp-Zeil}
 J.~Wimp and D.~Zeilberger,  Resurrecting the asymptotics of linear recurrences,  {J. Math. Anal.
Appl.} \textbf{111} (1985), 162--176.

\bibitem{xin-iterate}
G. Xin,  A fast algorithm for MacMahon's partition analysis,
 {Electron. J. Combin.}   {\textbf{11}} (2004) R58.

\bibitem{xin-residue}
G.~Xin,  A residue theorem for {M}alcev-{N}eumann series,
 {Adv. Appl. Math.}   \textbf{35} (2005), 271--293.


\bibitem{xin-weyl} G. Xin, Determinant Formulas Relating to Tableaux of Bounded Height,  {Adv.
Appl. Math.}  to appear, {\tt arXiv:0704.3381}.


\bibitem{creat-telescoping} D. Zeilberger, The method of creative telescoping,  {J. Symbolic
Comput.}   \textbf{11} (1991), 195--204.


\end{thebibliography}
\end{document}